\documentclass{amsart}
\usepackage{amstext,amsfonts,amssymb,amscd,amsthm,a4}
\usepackage{enumerate}

\newtheorem{theo}{Theorem}[section]
\newtheorem{lem}[theo]{Lemma}
\newtheorem{rem}{Remark}[section]
\newtheorem{defn}{Definition}[section]
\newtheorem{prop}{Proposition}[section]
\newtheorem{cor}[theo]{Corollary}
\newtheorem{ex}{Example}[section]
\newtheorem{proper}{Properties}

\newcommand{\pl}{\mathbb{P}}

\newcommand{\cl}{\mathbb{C}}

\newcommand{\pq}{\nabla}

\newcommand{\lrarrow}{\longrightarrow}
\newcommand{\fs}{\frac{dy-p\;dx}{\parp(F)}}
\newcommand{\fsk}{\frac{dy-p\;dx}{\parp(F_k)}}

\newcommand{\X}[1]{\dx(#1)+p\,\dy(#1)}

\newcommand{\Mr}{\mathcal{M}}
\newcommand{\deuf   }{dx\wedge dy}

\newcommand{\la }{\left\{ }
\newcommand{\ra }{\right\} }
\newcommand{ \rp }{\right) }
\newcommand{ \lp }{\left( }

\newcommand{\ol}{\mathcal{O}}
\newcommand{\dx}{\partial_x}
\newcommand{\dy}{\partial_y}

\newcommand{\parp}{\partial_p}
\newcommand{\W}{\mathcal{W}}

\newcommand{\dept}{~\upshape:\itshape}
\newcommand{\pvir}{~\upshape;\itshape}
\newcommand{\cro}[1]{\left[#1\right]}
\newcommand{\pare}[1]{\left(#1\right)}
\newcommand{\ak}[1]{\la #1\ra}

\def\mapdown{\Big\downarrow}

\begin{document}
\title[The connection associated with webs]{Properties of the connection associated with planar webs and applications}
\author{Olivier Ripoll}
\curraddr{Universit\'e Bordeaux I, Institut de Math\'ematique
Bordelais} \email{Olivier.Ripoll@math.u-bordeaux1.fr}
\address{I.M.B., 351 cours de la Lib\'eration, 33405 Talence}

\subjclass[2000]{53A60 ; 34A20 ; 58A15}

\date{April 15, 2007}

\keywords{Web geometry, differential equations, connections,
Cartan-Spencer theory}

\begin{abstract}
We give various results and applications using the connection
$(E,\nabla)$ associated with a~$d$-web in \cite{H-1}. More
precisely, we exhibit fundamental $\ol^*$-invariants of the web
related to the differential equation of first order which presents
the web. They cast some new lights on the connection and its
construction, both conceptually and effectively. We describe~$4$
and~$5$-webs from this point of view and show for instance that
the connection gives account for the linearizability conditions of
the web. Moreover, we get characterization of maximal rank webs
such as exceptional~$5$-webs, and $4$-webs \emph{via} a new proof
of the Poincar\'e theorem in terms of differential systems. We
establish the trace formula related to the determinant bundle
$(\det E, \det \nabla)$ and the extracted $3$-webs. Furthermore,
the theorem of determination of the rank is proved to give an
explicit criterion for measuring the rank of a web.
\end{abstract}
\maketitle
\section{Introduction}

Planar webs geometry is dedicated to the study of classes of first
order ordinary differential equations  $$F(x,y,y'):=a_0(x,y)\cdot
(y')^d+a_1(x,y)\cdot (y')^{d-1}+\cdots+a_d(x,y)=0$$ with
coefficients in the ring~$\ol=\cl\ak{x,y}$ of convergent power
series in two variables, up to an element in the group~$\ol^*$ of
invertible in $\ol$. Particulary, we are interested in the study
of specific relations between the solutions of these equations,
called the abelian relations.

These differential equations of degree~$d$ admit $d$ solutions
out of the singular locus given by their $y'$-resultant. They
are called the leaves of the web.

A basic problem in web geometry is to count the number of linearly independent abelian relations, which is the rank of the web, and to determine the webs of maximal rank. The connection $(E,\pq)$ we deal with here answers to this last problem, but not only, as we will prove it in this paper. Questions about the simultaneous linearizability of the leaves of the web also appear naturally and we will show how they relate to the first ones.

According to the classical definition of webs, the leaves of a planar~$d$-web
can equivalently be given by $d$ foliations defined by the level
sets~$F_i(x,y)=constant$ for~$1\leq i\leq d$ where
the~$F_i$ are in~$\ol$ in general position, with $F_i(0)=0$.
An abelian relation is then a relation of the
form~$$\sum_{i=1}^dg_i(F_i)dF_i=0$$ where the~$g_i$ are analytic
in one variable.

 $ $

The pioneers of web geometry are W. Blaschke and G. Bol in the thirties
(\cite{BB}). Then the subject get another raise thanks to S. S. Chern
and P. A. Griffiths in the seventies (\cite{CG}). Nowadays, the
subject is still up to date with the works of M. A. Akivis,
V. V. Goldberg and V. V. Lychagin (See for instance \cite{GL}), I.
Nakai, E. Ghys, D. Cerveau and A. H\'enaut, more recently D.
Lehmann and V. Cavalier \cite{CL}, L. Pirio and J.-M. Tr\'epreau
(\cite{P} or \cite{PT}), G. Robert, J. V. Pereira and D. Marin
\cite{PPM}.

Thus, the subject has drawn links with many subjects such as
foliations,~$\mathcal{D}$-modules, polylogarithms, differential
equations and Cartan-Spencer techniques for instance and could
meet now some topics like differential algebra (see for instance
\cite{RS} on singular solutions of planar webs or \cite{RS2} on
Darboux polynomials), differential Galois theory and meromorphic
connections (see \cite{H-2}). But web geometry has also some
applications in applied mathematics. Non linear optic geometry
(see for instance the article of J.-L. Joly, G. M\'etivier and J.
Rauch in \cite{W}) but also economy use web geometry where the
question of finding the rank of the web is crucial.

$ $

This paper aims to give new conceptual and effective results by
using the connection associated with webs. Let us describe the content of this paper. We first introduce shortly the tools used here by giving a description of the abelian relations of a web \emph{via} the data of the first order differential equation which presents the web. The connection is also described.

We then introduce the associated polynomials in theorem \ref{uv}
which lead us to the fundamental $\ol^*$-invariants of the web.
They are closely linked with the geometry of the differential
equation. Those $\ol^*$-invariants will give significative results
in association with the connection. But first, we show how they
relate to linear and algebraic webs and how they can be linked
with abelian relations. We then reduce their number by showing
that they can only be given by the data of a particular polynomial
and a fundamental $1$-form. These will be at the centre of this
study.

We then fully construct the connection associated with~$4$ and~$5$-webs. We show that it gives
account of the linearizability conditions (theorem \ref{linea}) and makes it clearer that if the Poincar\'e theorem is true for $4$-webs, as we prove it in theorem \ref{poin}, it is
not the same for $5$-webs. By the way, the obstruction is given in
this case (theorem \ref{vex}). Hence, we give a new
conceptual description of exceptional~$5$-webs (proposition
\ref{5w}).

The trace formula is then proved (theorem \ref{trace}). For this, we need a
few propositions which state the link between the connections
associated with extracted webs and the connection of the initial
web. They focus on the nature of the abelian relations. Then the proof is completed.

We end with the theorem of determination of the rank (theorem
\ref{rk}). It gives explicitly the locally free~$\ol$-modules of
finite type deduced from the local system of the abelian
relations. This fibre bundle is thus generated by the abelian
relations of the web. The particular form of the connection allows
to embody this $\ol$-modules so that it is now possible to compute
the rank without searching the abelian relations, as it was done
before in general. Then an exemple is given.

\section{Objects and tools}\label{ot}
References for this introduction of our objects and tools are the
founder books \cite{BB} or \cite{B} of W. Blaschke and G. Bol, the
article of S. S. Chern \cite{C} and the book \cite{W}.\ Unless
further specifications, the notations adopted for the sequel of this
article are fixed in this part.

\begin{defn}A non singular planar~$d$-web~$\W(d)$ is defined by the family of the
solutions of an ordinary differential equation of the first order
$$F(x,y,y'):=a_0(x,y)\cdot (y')^d+a_1(x,y)\cdot
(y')^{d-1}+\cdots+a_d(x,y)=0$$ which is a~$y'$-polynomial with
coefficients in~$\ol$, in a neighborhood of a point where the
$y'$-resultant~$R(x,y)$ of~$F$ and~$\partial_{y'}(F)$ is non zero.
\end{defn}
We say that a~$d$-web is \textit{presented} by such a differential
equation. In the sequel, a $d$-web will always be presented by
such a differential equations, denoted $F=0$, if it is not
otherwise specified. An $\ol^*$-invariant of a web~$\W(d)$ only
depends on the class of differential equations, modulo an
invertible in $\ol$ and not on a particular equation in this
class. We could easily justify this definition of invariance by
saying that the natural condition we could expect from our objects
is that as~$F=0$ and~$\rho.F=0$ present the same web if~$\rho$ is
invertible in~$\ol$, they do not depend on the presentation of the
web.

A theorem of Cauchy asserts that in a neighborhood of a point
$(x_0,y_0)$ such that~$R(x_0,y_0)\neq0$, the equation admits~$d$
integral curves, which are the leaves of the~$d$-web. Close
to this point, the polynomial~$F(x,y,p)$ admits~$d$ distinct roots
which are the slopes of the leaves, usually noted~$p_i(x,y)$ for
$1\leq i\leq d$. We will denote by~$\W(F_1,F_2,\ldots,F_d)$
a~$d$-web given by the level functions~$F_1,F_2,\ldots,F_d$, as it
was defined in the introduction. Thus, if we set $p_i:=-\dx(F_i)/\dy(F_i)$, then $\W(F_1,F_2,\ldots,F_d)$ is presented by the differential equation $\prod_{i=1}^d(y'-p_i)=0$.

We will use here the word web instead of non singular web, which is more convenient, since we will always be in the neighborhood of a point where the web is defined, that is to say where the $y'$-resultant of $F$ is non zero.
\begin{defn}The $\cl$-vector space defined by
$$\mathcal{A}(d)=\big\{ \pare{g_1(F_1),\ldots,g_d(F_d)}\in
\ol^d\mbox{ with }g_i\in\cl\la t\ra \mbox{ and
}\sum_{i=1}^dg_i(F_i)dF_i=0\big\}$$is called the space of
\textit{abelian relations} of the~$d$-web~$\W(d)$ \end{defn} We
have the following theorem:
\begin{theo}[Bol 1932, Blaschke 1933]The dimension of~$\mathcal{A}(d)$ is finite and called the
\textit{rank of the web}. It is an $\ol^*$-invariant of the web
and more, we have the following optimal inequality:
$$\mbox{\upshape rk \itshape}\W(d):=\dim_\cl \mathcal{A}(d) \leq \frac{1}{2}\,(d-1)(d-2).$$
\end{theo}
Note that in the sequel, $d$ will always be an integer greater or
equal to $3$, since abelian relations are trivial otherwise. We
denote by $\pi_d$ the integer $\frac{1}{2}\,(d-1)(d-2)$.

Taking an example, consider the~$3$-web~$\mathcal{H}=\W(x,y,x+y)$, given by
the level sets~$x=constant, y=constant\mbox{ and }x+y=constant$.
The rank of $\mathcal{H}$ is maximal, equal to $1$
since~$dx+dy-d(x+y)=0$ is a non trivial abelian relation.

$~$

For a~$3$-web~$\W(F_1,F_2,F_3)$, we know with the work of Blaschke
that there exists a differential~$1$-form~$\gamma$ with
coefficients in $\ol$ such that $d\omega_i=\gamma\wedge\omega_i$
where~$\omega_i=\rho_idF_i$, with~$\rho_i\in\ol^*$ chosen such
that the normalization~$\omega_1+\omega_2+\omega_3=0$ holds.
The~$2$-form~$d\gamma$ is, unlike~$\gamma$, an $\ol^*$-invariant
of the web called the Blaschke curvature of~$\W(3)$. We then have
the following equivalence~$\mbox{rk }\W(3)=1$ \textit{if and only
if}~$d\gamma=0$.

Our interest is specially turned on natural classes of webs. For
instance, one can consider the class of \textit{linear webs}, whose
leaves are germs of straight lines. Linearizable webs are those
for which there exists a change of coordinates which transforms the
web in a linear one. Among them, one can consider algebraic webs.
Given a reduced algebraic curve~$C$ in~$\pl^2$ of degree~$d$, a
straight line cuts generically~$C$ in~$d$ points. By duality, we get a~$d$-web whose leaves are tangent to
the dual curve of~$C$, provided that~$C$ does not contain straight
lines. Webs constructed thanks to an algebraic curve are
\textit{algebraic} ones and algebrizable webs are those for which
there exists a change of coordinates which transform them in
algebraic ones. One can see that Abel's theorem states that the
rank of algebraic webs is maximal. Thanks to an
\textit{Abel-inverse} type theorem (\emph{cf.} \cite{G}), we get
the following:
\begin{theo}[Lie-Darboux-Griffiths]\label{LDB}
A linear~$d$-web~$\mathcal{L}(d)$ admitting an abelian relation
whose terms are all different from zero, is algebraic. In
particular, a linearizable~$d$-web with maximal rank is
algebrizable.
\end{theo}
A web of maximal rank is not always linearizable
if~$d$ is greater or equal to~$5$. Such a web will be called
\textit{exceptional}. Several exceptional webs have been
discovered in the last years by Robert and Pirio-Tr\'epreau for instance (see
\cite{P} and \cite{PT}), with the consequence to tone down the
range of the word exceptional but among them, one remains so. It is the first
example of such a web, given by G. Bol in 1936. The leaves of this~$5$-web
at a point~$z$ are given by~$4$ pencils of lines in
general position and the only conic passing through the~$4$ points
and~$z$. This web is associated with the functional relation of
the dilogarithm. It is an exceptional web and among them, it is the only one
to be also hexagonal: all its extracted 3-webs are of maximal
rank. One can refer to the article of P. A. Griffiths
\cite{G2} for a prospective survey of the link between
polylogarithms and exceptional webs.

With the setting adopted here, the $\cl$-vector
space $\mathcal{A}(d)$ of the abelian relations of a web is described
in \cite{H-1} as follow: let~$S=\ak{F(x,y,p)=0}$ be the surface
defined in~$\cl^3$ by $F$ with the projection~$\pi$ on the $(x,y)$
plane. There exists a~$\cl$-isomorphism between~$\mathcal{A}(d)$
and the space of $1$-forms~$$\mathfrak{a}_F=\big\{\omega=(b_3\cdot
p^{d-3}+b_4\cdot
p^{d-4}+\cdots+b_d)\cdot\displaystyle{\frac{dy-pdx}{\parp(F)}}\in\pi_\ast(\Omega_S^1),\;b_i\,\in\,\ol\mbox{
and }d\omega=0\big\}$$ where~$(\Omega^\bullet_S,d)$ is the
usual de Rham complex, with
$\Omega^\bullet_S=\Omega^\bullet_{\cl^3}/\pare{dF\wedge\Omega^{\bullet-1}_{\cl^3},F.\Omega^\bullet_{\cl^3}}$

Such a~$1$-form~$\omega=r.\frac{dy-pdx}{\parp(F)}$ is
in~$\Omega^1_S$ if and only if there exists elements $r_p$
in $\ol\cro{p}$ of degree less or equal to~$d-1$ in $p$, and
$t:=t_2\cdot p^{d-2}+\ldots+t_d$ with coefficients in $\ol$ such that
$$r.(\dx(F)+p.\dy(F))+r_p.\parp(F)=(\dx(r)+p.\dy(r)+\parp(r_p)-t).F$$
One then gets with the previous relation that
$$d\omega=t\cdot\frac{dy-pdx}{\parp(F)}\cdot$$
So, the form is closed if and only if~$t=0$ that is to say that
there exists $r_p$ in $\ol\cro{p}$ of degree at most $d-1$ such
that
$$\pare{\star} \qquad r.(\dx(F)+p.\dy(F))+r_p.\parp(F)=(\dx(r)+p.\dy(r)+\parp(r_p)).F$$
holds. It is also the same to say that the~$b_i$ are analytic
solutions of the following homogeneous differential systems~$\mathcal{M}(d)$, where the left members are the coefficients $t_2,\ldots,t_d$ of the polynomial $t$, equal to zero in this case :
$$\mathcal{M}(d)\;\;\left\{\begin{matrix}\dx (b_d) &+&  A_{1,1}\cdot b_3 + \cdots + A_{1,d-2}\cdot b_d&=&0\cr
\dx (b_{d-1})+\dy(b_{d})&+&   A_{2,1}\cdot b_3 + \cdots +
A_{2,d-2}\cdot b_d&=&0\cr
  & & & \vdots& \cr
\dx (b_{3})+\dy(b_{4})&+&  A_{d-2,1}\cdot b_3 + \cdots +
A_{d-2,d-2}\cdot b_d&=&0\cr \dy(b_3)&+&A_{d-1,1}\cdot b_3 +
\cdots + A_{d-1,d-2}\cdot b_d&=&0\end{matrix}\right.$$ The coefficients~$A_{i,j}$ are in~$\ol\cro{1/\Delta}$, where~$\Delta$
is the~$y'$-discriminant of~$F$. So the space of abelian relations
of the web is identified with the space of solutions of the
differential system~$\Mr(d)$. By the nature of this system, it is
actually a local system, and so for the space of abelian
relations, seen as a sheaf over $\cl^2-\ak{\Delta=0}$.

$~$

The exterior differential on~$S$ induces a linear differential
operator
$$\begin{matrix}\rho
:&\ol^{d-2}&\longrightarrow&\ol^{d-1}\cr&\pare{b_3,\ldots,b_d}&\longmapsto&\pare{t_2,\ldots,t_d},\end{matrix}$$
the~$t_i$ being given by the~$d-1$ equations of~$\mathcal{M}(d)$.

The study of this system leads us to consider, with the notations
defined in \cite{BC3G}, the jet space~$J_k(\ol^{d-2})$ of
order~$k$ over $\ol^{d-2}$
and~$j_k:\ol^{d-2}{\rightarrow}J_k(\ol^{d-2})$ the natural
derivation map. From this differential operator, one gets
the~$\ol$-morphism~$p_0:J_1(\ol^{d-2}){\rightarrow}\ol^{d-1}$ and
its successive prolongations~$p_k$. Let~$R_k$ be the kernel of
$p_k$. The Spencer complex associated with the prolongations~$p_k$
is then given by~$$0\rightarrow \mbox{Sol
}\mathcal{M}(d)\stackrel{j_{k+1}}{\rightarrow}R_k\stackrel{D}{\rightarrow}\Omega^1\otimes_\ol
R_{k-1}\stackrel{D}{\rightarrow}\Omega^2\otimes_\ol
R_{k-2}{\rightarrow}0.$$ Using properties of these objects which
will be detailed in the section~$3$ for~$4$-webs, the main result
of \cite{H-1} proves the existence of a~$\cl$-vector
bundle~$E:=R_{d-3}$ included in $J_{d-2}(\ol^{d-2})$ of
rank~$\pi_d$ on~$(\cl^2,0)$ which admits a connection
$$\nabla:E\longrightarrow\Omega^1\otimes_\ol E$$ such that the
space~$\mbox{Ker } \nabla$ of its horizontal sections is
isomorphic to $\mathcal{A}(d)$. Moreover, the curvature of $(E,\pq)$ takes its
values in $\Omega^2\otimes_\ol g\subset\Omega^2\otimes_\ol E$
where~$g$ is a free $\ol$-module of rank one over~$(\cl^2,0)$.

Note that this connection is a meromorphic connection, with poles
on the $y'$-discriminant of $F$.

In the case of $3$-webs, we get a natural normalization of the
webs by considering special $1$-forms on the surface $S$, so that
we can prove that the curvature of the connection thus constructed
is exactly the Blaschke curvature of the $3$-web. So the curvature
extends for $d$-webs the Blaschke curvature of $3$-webs. In the
sequel, referring to the Blaschke curvature of a $3$-web will mean
that we consider the curvature of the connection associated to a
$3$-web.

\section{Associated polynomials}
\subsection{Introduction}

We are looking for $\ol^*$-invariants both linked to the
differential equation and the abelian relations of the web
presented by this equation. For instance, this will potentially
allows us to join together conditions of ranking and
linearizability. Still, we want to get a minimal system of
$\ol^*$-invariants which will describe the web as completely as
possible. These $\ol^*$-invariants will be deduced from the so
called associated polynomials, whose existence is given by the
following theorem.
\begin{theo}[Associated polynomials]\label{uv}
Let~$\W(d)$ be a~$d$-web presented by a differential equation $F=0.$
There exists two polynomials~$U$ and~$V$ in $p$ respectively of
degree~$d-2$ and~$d-1$ whose coefficients admit poles on
$R=\mbox{Result}(F,\partial_{y'}(F))$, such that the following
equalities hold\dept
$$(\diamond)\qquad\partial_x(F)+p \partial_y(F)={U}\cdot F+{V}\cdot \partial_p(F)$$
Moreover such an expression is unique : if two
polynomials~$\widetilde{U}$ and~$\widetilde{V}$ in $p$
respectively of degree~$d-2$ and~$d-1$ satisfy $(\diamond)$, then
$U=\widetilde{U}$ and $V=\widetilde{V}$.\end{theo}

\begin{proof}Such polynomials~${U}:={u_2}\cdot p^{d-2}+\ldots+{u_d}$ and~$V:=v_1\cdot p^{d-1}+\ldots+{v_d}$
must satisfy a system~$S(\diamond)$, deduced from $(\diamond)$.
If~$\mathcal{R}$ stands for the Sylvester square matrix of
order~$2d-1$ whose determinant is $R$, the system is given by
$$S(\diamond)\qquad\mathcal{R}\cdot ^t\left(\begin{matrix}u_2,&
\ldots & u_d, & v_1, & \ldots & v_d\end{matrix}\right)$$
$$=^t\left(\begin{matrix}0, & \ldots & 0, &\dy(a_0),&
\dx(a_0)+\dy(a_1), & \ldots & \dx(a_{i-1})+\dy(a_i), & \ldots
&\dx(a_d)\end{matrix}\right).$$ We then get~$(\diamond)$ and its
uniqueness \emph{via} Cramer's rule.
\end{proof}
In this article we fix the following notations:$$U={u_2}\cdot
p^{d-2}+\ldots+{u_d}\mbox{  and }V:=v_1\cdot
p^{d-1}+\ldots+{v_d}.$$
\begin{cor}For all~$1\leq i\leq d-3$,
there exists two polynomials~$U_i$ and~$V_i$ respectively of
degree~$d-2$ and~$d-1$ whose coefficients admit poles on $R$, such
that the following equality hold\dept
$$(\diamond_i)\qquad p^i\cdot(\partial_x(F)+p \partial_y(F))={U_i}\cdot F+{V_i}\cdot \partial_p(F)$$
Moreover such an expression is unique.\end{cor}

The proof of this corollary is the same as the previous one.
For~$0\leq i\leq d-3$, the couples of
polynomials~$${U_i}:={u^i_2}\cdot p^{d-2}+\ldots+{u^i_d}\mbox{ and
}V_i:=v^i_1\cdot p^{d-1}+\ldots+{v^i_d}$$ are called the
associated polynomials of~$F$ of order~$i$, the associated polynomials of order $0$ being the couple $(U,V)$.

\subsection{Associated polynomials as web's $\ol^*$-invariants}
The associated polynomials are linked to the web thanks to the
following properties:
\begin{prop}\label{inv} Let~$\W(d)$ be a~$d$-web presented by
$F=0$ and let ~$(U_i^{F},V_i^{F})$ be the associated polynomials
of order~$i$. If~$(U_i^{g\cdot F},V_i^{g\cdot F})$ are the
associated polynomial of the equation~$g.F=0$, where~$g$ is
invertible in~$\ol$, we then have the following relations\dept
$$\left\{
\begin{matrix} U_i^{g\cdot F}&=&U_i^{F}+\frac{1}{g}\cdot
p^i\cdot(\partial_x(g)+p\partial_y(g))\cr \left.\right.\cr
V_i^{g\cdot F}&=&V_i^{F}\end{matrix} \right.$$\end{prop}
\begin{proof}
We need to write the relation we get from theorem \ref{uv}:\\
$p^i(\dx(g\cdot F)+p\dy(g\cdot F))=U_i^{g\cdot F}\cdot g\cdot
F+V_i^{g\cdot F}\cdot \partial_p(g\cdot F)$,
that is to say \\  $p^i(\dx(g)+p\dy(g))F+p^i(\dx(F)+p\dy(F))=g U_i^{g\cdot F}\cdot
F+V_i^{g\cdot F}\cdot g\partial_p(F)$. The equalities result from the uniqueness obtained in
theorem \ref{uv} and its corollary.
\end{proof}
So, the polynomials~$V_i$ are $\ol^*$-invariants of the web, as
the following coefficients of~$U_i$:
$u^i_2,\ldots,\widehat{u^i_{d-i-1}},\widehat{u^i_{d-i}},\ldots,u^i_{d}$
where the hat means that the coefficients are omitted.

Moreover, for~$0\leq k\leq  d-2$, the~$(d-3)(d-2)$ differences
$$\left\{\begin{matrix}u^{i-1}_{d-i+1}-u_{d}\qquad \mbox{for }
2\leq i\leq  d-2 \cr \left.\right. \cr
u^{j-2}_{d-j+1}-u_{d-1}\qquad \mbox{for } 3\leq j\leq
d-1\end{matrix}\right.~$$ and the~$d-2$ forms
$d(u^i_{d-i}dx+u^i_{d-i-1}dy)$ are $\ol^*$-invariants of the web.

\subsection{Linear and algebraic webs}
We will begin this section with results concerning the
linearizability of the webs.

Let~$\W(d)$ be a planar~$d$-web whose slopes of
the leaves are denoted for~$1\leq i\leq d$ by~$p_i\in\ol$ . There exists a
unique polynomial of degree less or equal to~$d-1$
$$P_{\W(d)}:=l_1\cdot p^{d-1}+l_2\cdot p^{d-2}+\cdots+l_d$$ whose
coefficients are in~$\ol$, such that the following equality holds
for all~$i$:
$$X_i(p_i):=\dx(p_i)+p_i\dy(p_i)=P_{\W(d)}(x,y,p_i(x,y)).$$

Then, the graphs of the leaves are solutions of the
equation~$y''=P_{\W(d)}(x,y,y')$. The following properties of
$P_{\W(d)}$ are due to H\'enaut:
\begin{proper}\label{pw}\normalfont
\begin{enumerate}
    \item The web~$\W(d)$ is linear if and only if $X_i(p_i)=0$ for
    all~$i$, if and only if ~$P_{\W(d)}=0$\pvir \\~$~$

    \item For~$d\geq 4$, the web~$\W(d)$ is linearizable if and only
if~$deg(P_{\W(d)})\leq 3$ \\ and~$(l_d,l_{d-1},l_{d-2},l_{d-3})$
is a solution of the non linear differential system\dept
\end{enumerate}     \vspace{1mm}
$\la\begin{matrix}L_1=-\dx(\dx(l_{d-2})-2\dy(l_{d-1}))-l_{d-1}(\dx(l_{d-2})-2\dy(l_{d-1}))-3\partial^2_y(l_d)
\cr  -3\dy(l_{d-2}l_d)+3\dx(l_dl_{d-3})+3l_d\dx(l_{d-3})=0\qquad\qquad\cr\cr
L_2=\dy(2\dx(l_{d-2})-\dy(l_{d-1}))-l_{d-2}(2\dx(l_{d-2})-\dy(l_{d-1}))-3\partial^2_x(l_{d-3})
\cr  +3\dx(l_{d-1}l_{d-3})
-3\dy(l_dl_{d-3})-3l_{d-3}\dy(l_d)=0.\qquad\qquad\end{matrix}\right.$
\end{proper}
The polynomial~$P_{\W(d)}$ is hence called the \emph{linearization
polynomial}. It appears to be a important $\ol^*$-invariant of the
web, as we will see. One can compare this property with webs
linearization results of Akivis, Goldberg and Lychagin (\emph{cf}.
 \cite{Go}).

\begin{prop}\label{pwv}
For~$0\leq k\leq d-3$ and all~$1\leq i\leq d$, we have the
following equality:~$V_k(x,y,p_i)=-(p_i)^k.X_i(p_i)$. In
particular,~$P_{\W(d)}=-V$
\end{prop}
\begin{proof}Since~$F=\prod_i^d(p-p_i)$ presents the web, we have
the following equalities:\\~$p^k.(\X F)=p^k.\sum_{i=1}^d-(\X
{p_i})\prod_{j=1,j\neq i}^d(p-p_j)$ and \\ $p^k.(\X
F)=U_k.F+V_k\sum_{i=1}^d\prod_{j=1,j\neq i}^d(p-p_j).$ Let
$p=p_i$, we have then~$V_k(x,y,p_i)=-(p_i)^k.X_i(p_i)$. By the
uniqueness of~$P_{\W(d)}$, the second equality follows.
\end{proof}
As a consequence of the properties of linear
webs given in the beginning of this subsection, we have the following corollary:
\begin{cor}[Linear webs]\label{lin}
A~$d$-web is linear if and only there exists an integer $0\leq
k\leq d-3$ such that~$V_k=0$.
\end{cor}
These are geometric properties, since they are expressed with
$\ol^*$-invariants of the web $\W(d)$.

$ $

It can be checked that an algebraic web is
presented by an equation~$$F(x,y,y')=g\cdot P(y-y'x,y')=0$$ where~$P
\in \cl \left[s,t\right]$ is an affine equation of the reduced
algebraic curve which defines the web, and~$g$ is invertible in
$\ol$. As linear webs, the polynomials~$V_i$ are equal to zero,
and we get an additional condition for such a web, given in the following theorem:
\begin{prop}[Algebraic webs]\label{propA} Let~$\W(d)$ be a~$d$-web
presented by a differential equation $F(x,y,y')=0$ and let~$(U,V)$ be its associated polynomials. We have the following equivalences\dept
\begin{enumerate}
\item[i)]~$\W(d)$ is algebraic\pvir
\item[ii)] There exists~$\phi$ in~$\ol$  such that~$V=0$ and
$U=\dy(\phi)p+\dx(\phi)$\pvir
\end{enumerate}
\end{prop}
\begin{proof}If~$\W(d)$ is algebraic, an equation
~$F(x,y,p)=e^\phi P(y-px,p)$ presents the web, where~$\phi$ is in~$\ol$ and~$P$
in~$\cl\left[s,t\right]$. So~$\dx (P(y-px,p))+p\dy (P(y-px,p))=0$
and, according to theorem \ref{uv} we get the properties~$ii)$.

Conversely, let~$G(x,y,p)=e^{-\phi}F(x,y,p)$. Since~$V(F)=0$ by hypothesis, then\\$V(G)=0$. Moreover, we
have~$U(G)=U(F)+e^\phi (\X {e^{-\phi}})=0$ and
so~$\X G=0.$ The polynomials~$y-px$ et~$p$ are linearly independent solutions of the previous equation. By Frobenius theorem, there exists an analytic function
~$\gamma$ in~$\ol$ such that~$G(x,y,p)=\gamma(y-px,p).$
In fact,~$\gamma$ is a polynomial in~$p$ of degree $d$ since we can show that for
its partial derivatives of order greater than $d$ are equal to zero.
The second equivalence is a consequence of the preceding one.
\end{proof}
This proposition is closely related to the  Lie-Darboux-Griffiths
theorem. The condition that the web is linear is expressed by the
cancelation of~$V$. The maximal rank condition can not be yet
interpreted as the condition on~$U$, but we will see that it is
the case (\emph{cf}. proposition \ref{alglin}).
\begin{rem}\normalfont The associated polynomials $(U,V)$ allow us
to find specific singular solutions of the differential equation
$F=0$. Moreover, in the case where such a singular solution exist,
the connection does not admit poles on the locus of this solution.
These results were obtained by J. Sebag and the author in
\cite{RS} by crossing the classical results on singular solutions
given by G. Darboux, and those of J. F. Ritt and E. R. Kolchin in
differential algebra, around the $y'$-resultant of $F$ which is
omnipresent in our constructions.
\end{rem}
\subsection{Link with the abelian relations}
The following lemma will make a bridge between those associated
polynomials and the abelian relation of the web.

\begin{lem}\label{link} With the previous notations, let~$r:=r(x,y,p)=b_3\cdot
p^{d-3}+\ldots+b_d$ with coefficients in~$\ol$. The polynomials $U_r=b_3\cdot U_{d-3}+\ldots+b_d\cdot U$ and
$V_r=b_3\cdot V_{d-3}+\ldots+b_d\cdot V$ are such that
$r\cdot(\partial_x(F)+p
\partial_y(F))=U_r\cdot F+V_r\cdot
\partial_p(F).$\\ Then, a~$1$-form
$$\omega=(b_3\cdot p^{d-3}+b_4\cdot
p^{d-4}+\cdots+b_d)\cdot\displaystyle{\frac{dy-pdx}{\parp(F)}}\in\pi_\ast(\Omega_S^1)\mbox { belongs to }\mathfrak{a}_F$$
if and only if~$r$ satisfies the
equation $U_r+\partial_p(V_r)=\dx(r)+p\dy(r)$.\end{lem}
Indeed,  the first
part is a direct consequence of the uniqueness property in theorem
\ref{uv} and its corollary. The relation $(\star)$ given in
section 2 gives the second part, again using the uniqueness
property.

This allows us to express the system~$\Mr(d)$ in terms of our
associated polynomials.

After a short computation we get the following expression for
$\Mr(d)$:
$$\left\{\begin{matrix}\dx (b_d) -  (u^{d-3}_d+v_{d-1}^{d-3})b_3 - \ldots - (u_d+v_{d-1})b_d&=&0\cr
   &\vdots&  \cr
\dx (b_{d+1-i})+\dy(b_{d+2-i})- \ldots -(u_{d+1-i}^{d-2-j}+i\cdot
v^{d-2-j}_{d-i})b_j  \ldots&=&0\cr  &\vdots&   \cr \dy(b_3)-
(u^{d-3}_2+(d-1)v_{1}^{d-3})b_3 - \ldots -
(u_2+(d-1)v_{1})b_d&=&0\cr\end{matrix}\right.$$

\begin{prop}Let $\alpha=A_1dx+A_2dy:=A_{1,d-2}dx+A_{2,d-2}dy$. We have the
following equality:~$$\alpha=(-\frac{\dx(a_0)}{a_0}-\dy(\frac{a_1}{a_0})+\sum_{i=1}^{d-1}v_i\sum_{k=1}^d
p_k^{d-1-i})dx+(-\frac{\dy(a_0)}{a_0}+\sum_{i=1}^{d-2}v_i\sum_{k=1}^d
p_k^{d-2-i})dy.$$
\end{prop}
The proof is a direct computation, using the writing of $\alpha$
with our $\ol^*$-invariants:\\
$\alpha=-(u_d+v_{d-1})dx-(u_{d-1}+(d-2)v_{d-2})dy$.

The form $\alpha$ will be called the \emph{fundamental}~$1$-form.
Its differential is an $\ol^*$-invariant of the web with
proposition \ref{inv}. This form will play an important part in
the sequel, justifying the distinction we make.

For~$d=3$, the system~$\Mr(3)$ is given by the
coefficients~$A_{11}=A_1$ and~$A_{21}=A_2$. The associated connection is a $1$-form, which, in a suitable basis is
the fundamental $1$-form~$\alpha=A_1dx+A_2dy$. It has been shown in
\cite{H-1} that its curvature is then the Blaschke
curvature~$d\gamma$ of the $3$-web, which emphasizes on the
importance of~$\alpha$. So, in the case of a $3$-web the
differential of the fundamental form is the Blaschke curvature.

$ $

We are now looking for a minimal set of $\ol^*$-invariants. The
next theorem gives us such a minimal set, which will play a
central part in this article.
\begin{theo}The system~$\Mr(d)$ and consequently, the connection, can only be written
thanks to~$\alpha$ and~$V$. Moreover, $d\alpha$ and $V$ are
$\ol^*$-invariants of the web.
\end{theo}
\begin{proof}  First we can deduce all the polynomials~$V_k$ from
the polynomial~$V$ and~$F$. Let us remark that for all~$1\leq
k\leq d-3$, we have the following polynomial equalities:
$$(V_k-pV_{k-1})(x,y,p)=-\frac{v^{k-1}_1}{a_0}F(x,y,p) \mbox{ and
}U_k-pU_{k-1}=\frac{v^{k-1}_1}{a_0}\parp(F)(x,y,p).$$ Indeed,
since the polynomials~$(V_k-pV_{k-1})(x,y,p)$ of degree~$d$ admits
the~$d$ slopes of the web~$p_i\in \ol$ as solutions, we get the first equality.
The second one is deduced from the uniqueness in theorem \ref{uv}.
So, if~$V=-P_{\W(d)}$ is known, we know all the others polynomials
$V_k$.

It remains to show that~$U$ can be deduced from~$V$. Considering
the system~$(\diamond)$ in theorem \ref{uv}, and taking~$V$ as a
parameter, the special form of the Sylvester determinant gives a
triangular system in the~$a_i$ and~$V$ with~$U$ as unknown. One
can check that using the Newton's relations between the
coefficients and the roots of an our equation $F$, the data
of~$\alpha$ and~$V$ allows us to compute all the coefficients of
the system~$\Mr(d)$.
\end{proof}

We will see in the next section the usefulness of this theorem,
and the special form of the coefficients of $\Mr(d)$ expressed
with $\alpha$ and $V$.  But first, to emphase on the fundamental
form, let us consider the linear case. We said that the
form~$\alpha$ is linked to the abelian relations. Precisely, we
have the following proposition which completes the result given in
proposition \ref{propA}:
\begin{prop}\label{alglin}The following properties are equivalent\dept
\begin{enumerate}
    \item[i)]$\W(d)$ is algebraic\pvir
    \item[ii)]$\W(d)$ is linear and the fundamental form satisfies $d\alpha=\partial_y^2(\frac{a_1}{a_0})=0$.
\end{enumerate}
\end{prop}
\begin{proof}
The coefficients of the system~$\Mr(d)$ for a linear~$d$-web is
the following, since the $V_k$ are zero:
$$\left(\begin{matrix}
0     & \ldots & 0      & 0        & -u_d    \cr 0     & \ldots &
0      & -u^1_{d-1}& -u_{d-1}\cr 0     & \ldots & \ddots &
-u^1_{d-2}& 0        \cr \vdots&        & \ddots &          & 0
\cr -u^{d-3}_3 &  0     &    0   & \ldots   & 0 \cr -u^{d-3}_2 & 0
& 0 &\ldots & \cr\end{matrix}\right).$$ But writing the~$u_j^i$ in
terms of the~$v_j^i$ shows that, since all the~$v_j^i=0$:
~$$u_d=u^1_{d-1}=\ldots=u^{d-3}_3=\frac{\dx(a_0)}{a_0}+\dy(\frac{a_1}{a_0})\mbox{
and }u_{d-1}=u^1_{d-2}=\ldots=u^{d-3}_2=\frac{\dy(a_0)}{a_0}.$$
\indent According to proposition \ref{propA} the linear web is
algebraic if and only if the following system admits a solution:
$$\la\begin{matrix}\dy(\phi)=\dy(a_0)/a_0\cr\dx(\phi)=\dx(a_0)/a_0+\dy(a_1/a_0)\end{matrix}\right.$$
\indent  The integrability condition is
then~$\partial_y^2(\frac{a_1}{a_0})=0$, that is to
say~$d\alpha=0$, which proves the proposition.
\end{proof}

\section{Geometric study of the~$4$-webs associated
connection}\label{co4} The general settings given in the preceding
sections will find a direct application in this part. We give here
an explicit computation of the connection associated with
a~$4$-web. Not only that it gives a computation tool, it will
allows us to find some more properties of the connection.
\subsection{Construction}
We adopt here the notations of the book \cite{BC3G} and the
methods developed in the article \cite{H-1}. In the case
of~$4$-webs, on can define the operator
$$\begin{matrix} p_0: &J_1(\ol^2)&\lrarrow& \ol^3\cr & \left(\begin{matrix}z_3,p_3,q_3 \cr
z_4,p_4,q_4\end{matrix}\right) &\lrarrow&
\left(\begin{matrix}p_4+A_{11}z_3+A_{12}z_4 \cr
p_3+q_4+A_{21}z_3+A_{22}z_4 \cr q_3+A_{31}z_3+A_{32}z_4
\end{matrix}\right)\end{matrix}.$$ \indent Then we get an exact and
commutative diagram

$$\begin{matrix}
   &  &  0\;\;\;\;\;\,&  &0\;\;\;\;\;\,  &   & 0 &   &   \cr
   &  &  \mapdown\;\;\;\;\;  &  & \mapdown\;\;\;\;\;  &   & \mapdown\;  &  & \cr
0  & \lrarrow & g_0\;\;\;\;\, &\lrarrow  & S_1(\ol^2)\;\,\,
&\stackrel{\sigma_0}\lrarrow   &  \ol^3  &\lrarrow   &  0 \cr
   &  &  \mapdown\;\;\;\;\;\,&  &  \mapdown\;\;\;\;\;\,&   &   \mapdown\;&   &   \cr
0  & \lrarrow & R_0\;\;\;\;\; &\lrarrow  & J_1(\ol^2)\;\,\,
&\stackrel{p_0}\lrarrow   & \ol^3  &\lrarrow   & 0  \cr
   &  &  \mapdown \overline{\pi}_{-1}&  &  \mapdown \pi_{-1}&   &   \mapdown\;&   &   \cr
0  & \lrarrow & R_{-1}\;\;\; &\lrarrow  & \ol^2\;\;\;\;
&\stackrel{p_{-1}}\lrarrow   & 0  &   &   \cr
   &  &  \mapdown\;\;\;\;\;\,&  &  \mapdown\;\;\;\;\;\,&   &   &   &   \cr
  &  & 0\;\;\;\;\;\, & & 0\;\;\;\;\;\, & &   & &   \cr
\end{matrix}$$ where the upper line concerns symbols of $p_0$.

Explicitly, we have: $\sigma_0(\begin{matrix}p_3,q_3\cr
p_4,q_4\end{matrix})=(p_4,p_3+q_4,q_3)$ whose kernel~$g_0$
is then isomorphic to $\ol$.

The first prolongation $p_1:J_2(\ol^2)\lrarrow J_1(\ol^3)$ of
$p_0$ is then given by
$$ p_1\left(\begin{matrix}z_3,p_3,q_3,r_3,s_3,t_3\cr z_4,p_4,q_4,r_4,s_4,t_4\end{matrix}\right)=
\lp\begin{matrix}p_4+A_{11}z_3+A_{12}z_4 \cr
p_3+q_4+A_{21}z_3+A_{22}z_4 \cr
q_3+A_{31}z_3+A_{32}z_4\end{matrix}\right.$$ $$
\left.\begin{matrix}r_4+A_{11}p_3+A_{12}p_4+\dx(A_{11})z_3+\dx(A_{12})z_4
\cr r_3+s_4+A_{21}p_3+A_{22}p_4+\dx(A_{21})z_3+\dx(A_{22})z_4 \cr
s_3+A_{31}p_3+A_{32}p_4+\dx(A_{31})z_3+\dx(A_{32})z_4\end{matrix}\right.$$
$$\left.\begin{matrix}s_4+A_{11}q_3+A_{12}q_4+\dy(A_{11})z_3+\dy(A_{12})z_4
\cr s_3+t_4+A_{21}q_3+A_{22}q_4+\dy(A_{21})z_3+\dy(A_{22})z_4 \cr
t_3+A_{31}q_3+A_{32}q_4+\dy(A_{31})z_3+\dy(A_{32})z_4
\end{matrix}\rp$$

which leads us to another commutative and exact diagram:
$$\begin{matrix}
   &  &  0\;\;\;&  &0\;\;\;  &   & 0  &   &   \cr
   &  &  \mapdown\;\;\;  &  & \mapdown\;\;\;  &   & \mapdown  &  & \cr
0  & \lrarrow & g_1=0 &\lrarrow  & S_2(\ol^2)
&\stackrel{\sigma_1}\lrarrow   &  S_1(\ol^3)  &\lrarrow   &  0 \cr
   &  &  \mapdown\;\;\;&  &  \mapdown\;\;\;&   &   \mapdown&   &   \cr
0  & \lrarrow & R_1\;\; &\lrarrow  & J_2(\ol^2)
&\stackrel{p_1}\lrarrow   & J_1(\ol^3)  &\lrarrow   & 0  \cr
   &  &  \mapdown \overline{\pi}_0&  &  \mapdown \pi_0&   &   \mapdown&   &   \cr
0  & \lrarrow & R_{0}\;\; &\lrarrow  & J_1(\ol^2)
&\stackrel{p_{0}}\lrarrow   & \ol^3  &   \lrarrow   & 0     \cr
   &  &  \mapdown\;\;\;&  &  \mapdown\;\;\;&   & \mapdown  &   &   \cr
  &  & 0 & & 0 & & 0  & &   \cr\end{matrix}.$$
Thus $\overline{\pi}_0$ is an isomorphism.
The first Spencer complex is given by~$$0\longrightarrow \mbox{Sol
} \mathcal{M}(4) \stackrel{j_2}{\longrightarrow} E
\stackrel{D}{\longrightarrow} \Omega^1\otimes_\ol R_0
\stackrel{D}{\longrightarrow} \Omega^2\otimes_\ol
R_{-1}\longrightarrow 0$$ where ~$E$ stands for the kernel of
$p_1$. We deduce from this the following commutative diagram, whose
lines are exact and whose columns are exact in~$R_k$, and~$j_2$
and~$j_3$ are monomorphisms:
$$\begin{matrix} & & &       & 0\;\;\;   & & 0\;\;\;   & & \cr
                    & & &               &    \mapdown\;\;\; & & \mapdown\;\;\; & & \cr
                    & & 0  &\longrightarrow & \mbox{Sol }\mathcal{M}(4)&\longrightarrow &\mbox{Sol }\mathcal{M}(4)&\longrightarrow & 0 \cr
                    & &  & & \mapdown j_3 & & \mapdown j_2 & &  \cr
                         & &    0  &\longrightarrow & R_2\;\; &\stackrel{\overline{\pi}_1}{\longrightarrow} & E=R_1 &\stackrel{\beta_1}\longrightarrow & \mathfrak{K}_{1} \cr
                        & & & & \mapdown\;\;\; & \swarrow_\pq & \mapdown D & &  \cr
                        & & 0  &\longrightarrow & \Omega^1\otimes_\ol E  &\stackrel{\overline{\pi}_0}\longrightarrow &\Omega^1\otimes_\ol R_0 &\longrightarrow & 0 \cr
                & &     & & \mapdown\;\;\; & & \mapdown D & &  \cr
                0   & \longrightarrow &  \Omega^2\otimes_\ol g_0  &\longrightarrow & \Omega^2\otimes_\ol R_0  &\stackrel{\overline{\pi}_{-1}}\longrightarrow &\Omega^2\otimes_\ol R_{-1} &\longrightarrow & 0 \cr
                    & &  & & \mapdown\;\;\; & & \mapdown\;\;\; & &  \cr
                     & & &       &    0\;\;\;   &    &   0\;\;\;   &     &      \cr\end{matrix}.$$ \indent
One can check that the kernel~$E$ of~$p_1$ is a $\cl$-vector bundle of
rank~$3$. The main
result of this construction is that, using properties of $D$, the
map~$\pq=\overline{\pi}_0^{-1}\circ D: E\lrarrow
\Omega^1\otimes_\ol E$ is a connection on $E$, and its kernel is
isomorphic to~$\mbox{Sol }\mathcal{M}(4)$.
\begin{rem}\normalfont Since~$\mbox{Ker }\pq$ is isomorphic to~$\mbox{Sol }\mathcal{M}(4)$ which is a local
system, the kernel of~$\mbox{Ker }\pq$ is also a local system,
even if the connection is not integrable.
\end{rem}
One can pick up an adapted basis of~$E=\mbox{Ker }(p_1)$. We
can choose such a basis $e_1,e_2,e_3$ this way, so that the
curvature takes its values in~$\Omega^2\otimes_\ol g_0$:
$$e_1=\left(\begin{matrix}0&-1&0&A_{21}+A_{12}&A_{31}&-A_{32}\cr
0&0&1&A_{11}&-A_{12}&-A_{22}-A_{31}\end{matrix}\right)$$

$$e_2=\left(\begin{matrix}-1&A_{21}&A_{31}&-\dy(A_{11})+\dx(A_{21})-A_{21}^2-A_{11}(A_{22}-A_{31})\cr
0&A_{11}&0&-A_{11}(A_{21}+A_{12})+\dx(A_{11})\qquad\qquad\qquad\qquad\end{matrix}\right.$$

$$\left.\begin{matrix}-A_{31}A_{21}-A_{32}A_{11}+\dx(A_{31})&-A_{31}^2+\dy(A_{31})\qquad\qquad\qquad\cr
-A_{11}A_{31}+\dy(A_{11})\qquad\qquad&-\dx(A_{31})+\dy(A_{21})+A_{32}A_{11}\end{matrix}\right)$$

$$e_3=\lp\begin{matrix}0&0&-A_{32}&\dy(A_{12})-\dx(A_{22})-A_{11}A_{32}\cr
 1&-A_{12}&-A_{22}&A_{12}^2-\dx(A_{12})\qquad\qquad\qquad
 \end{matrix}\right.$$
$$\left.\begin{matrix}A_{32}A_{12}-\dx(A_{32})\qquad\qquad&A_{32}(A_{31}+A_{22})-\dy(A_{32})\qquad\qquad\qquad\qquad\cr
A_{11}A_{32}+A_{12}A_{22}-\dy(A_{12})&A_{32}(A_{21}-A_{12})+A_{22}^2+\dx(A_{32})-\dy(A_{22})\end{matrix}\rp$$
   $ $

Since
$j_2(f)=\left(f,\dx f ,\dy f ,\dx^2 f ,\dx\dy f ,\dy^2 f \right),$
the exactness in~$E$ of the Spencer complex allows us to
compute~$D(e_1)$,~$D(e_2)$ and~$D(e_3)$ in $\Omega^1\otimes_\ol
R_0$, that must be composed with $\overline{\pi}_0$ to get the
value of~$\pq(e_i)$. Thus, we  get the connection matrix~$\gamma$
of~$\pq$ in this basis:

$$\gamma=\left(\begin{matrix} A_{12}dx+A_{31}dy & \xi_1 & \xi_2 \cr -dx &
A_{21}dx+A_{31}dy & -A_{32}dy \cr -dy & -A_{11}dx &
A_{12}dx+A_{22}dy\cr\end{matrix}\right)$$ where
$\xi_1 $ and $\xi_2$ are expression in the $A_{ij}$ which will be
detailed later. The matrix of the system~$\Mr(4)$ written with the
fundamental~$1$-form~$\alpha=A_1dx+A_2dy$ and~$V$ is
$$(A_{ij})= \quad \lp\begin{matrix}-v_4 & A_1 \cr A_1-v_3 & A_2 \cr A_2-v_2 & v_1\end{matrix}\rp.$$ Then we have the
following expression of the
connection matrix:$$\gamma=\left(\begin{matrix} A_1dx+(A_2-v_2)dy & \xi_1
& \xi_2 \cr -dx & (A_1-v_3)dx+(A_2-v_2)dy & -v_1dy \cr -dy & v_4dx
& A_1dx+A_2dy\cr\end{matrix}\right)$$ where
$\xi_1=(\dy(v_4)+v_4v_2)dx+(v_1v_4+\dx(A_2-v_2)-\dy(A_1-v_3))dy$
and\\ $\xi_2=(v_4v_1-(\dx(A_2)-\dy(A_1)))dx+(v_1v_3-\dx(v_1))dy$.

Its curvature is:$$d\gamma+\gamma\wedge\gamma=\left(\begin{matrix}
k_1 & k_2 & k_3 \cr 0 & 0 & 0 \cr0& 0 &
0\cr\end{matrix}\right)dx\wedge dy$$ where $k_1=d(\mbox{tr
}\gamma)$, $k_2$ and $k_3$ are computable but will be omitted
here (\emph{cf}. \cite{R2} for details).

\subsection{Interpretation}
We will give several applications of the writing of the connection
in terms of the~$1$-form~$\alpha$ and the polynomial~$V$. Since we have the equality~$V=-P_{\W(d)}$ established
in proposition \ref{pwv}, the linearizability conditions for
the web will be seen in the connection as we can see in the
following theorem.

\begin{theo}\label{linea}Up to a change of basis, the curvature matrix
associated with a~$4$-web is:
$$K=\left(\begin{matrix} k_1 & \dx(k_1)+L_1 & \dy(k_1)+L_2 \cr0 & 0 & 0 \cr0& 0 & 0\cr\end{matrix}\right)dx\wedge
dy$$ where~$L_1$ and~$L_2$ are  defined in properties \ref{pw}.
Thus, the curvature is an $\ol^*$-invariant of the web.
\end{theo}
\begin{proof}Using the identification~$V=-P_{\W(d)}$, and the expression of the curvature in terms of
$V$, one can compute that~$$k_2=\frac{1}{3}(\dx(k_1)+v_3k_1+L_1)
\mbox{ and }k_3=\frac{1}{3}(\dy(k_1)-v_2k_1+L_2).$$ The change of basis defined by the matrix
$$P=\lp\begin{matrix}\frac{1}{3} & -\frac{1}{3}v_3 &\frac{1}{3}v_2
\cr 0 & 1 &0\cr 0 &0 &1\end{matrix}\rp$$ gives the
needed expression.
\end{proof}
Now we can give  a new proof of a Poincar\'e theorem,
using only differential systems tools:
\begin{theo}[Poincar\'e, 1901]\label{poin}
A~$4$-web of maximal rank is linearizable.
\end{theo}
\begin{proof}Since the web is of maximal rank, the curvature is zero, thus,
~$k_1=k_2=k_3=0$. Hence~$L_1=L_2=0$ and since the degree
of~$P_{\W(4)}$ is~$3$, the web is linearizable by the properties
\ref{pw}.
\end{proof}
\section{The case of~$5$-webs}
In the case of~$5$-webs (and higher), computations are more
complicated, but the method is still the same and we can get
similar results than in the case of $4$-webs.
\subsection{The connection associated with~$5$-webs}
Let~$\W(5)$ be a planar~$5$-web, presented by a differential
equation $F=0$ of degree $5$. The matrix~$(A_{ij})$ of the
system~$\Mr(5)$ is
$$(A_{ij})=\lp\begin{matrix}\frac{a_5}{a_0}v_1 & -v_5 & A_1 \cr \cr
-2v_5+\frac{a_4}{a_0}v_1 & A_1-v_4&A_2\cr \cr
A_1-2v_4+\frac{a_3}{a_0}v_1 & A_2-v_3 &
2v_2-\frac{a_1}{a_0}v_1\cr\cr A_2-2v_3+\frac{a_2}{a_0}v_1 &
v_2-\frac{a_1}{a_0}v_1 & v_1\end{matrix}\rp$$ where we have
written $P_{\W(5)}=-v_1p^4-v_2p^3-v_3p^2-v_4p-v_5$ and the
fundamental form~$\alpha=A_1dx+A_2dy$. The construction of the
connection associated with the web gives us a $\cl$-vector bundle~$(E,\pq)$ of
rank~$6$ and a choice of an adapted basis of~$E$.

We can write the trace~$k_1$ of the curvature matrix :
$$k_1=6(\dx(A_2)-\dy(A_1))+4\dy(v_4)-8\dx(v_3)+3\dx(v_1\frac{a_2}{a_0})-\dy(v_1\frac{a_3}{a_0})$$
which is still an $\ol^*$-invariant of the web. With the notations
of the properties \ref{pw}, we have by computation the  following
theorem, the analogous  of theorem \ref{linea} for $5$-webs:
\begin{theo} The curvature is given in a suitable basis by the
matrix:
\small$$K=\left(\begin{matrix}k_1&\widetilde{k_2}+<v_1>_2&\widetilde{k_3}+<v_1>_3&\widetilde{k_4}+<v_1>_4&\widetilde{k_5}+<v_1>_5&\widetilde{k_6}+<v_1>_6\cr0&0&0&0&0&0\cr0&0&0&0&0&0\cr0&0&0&0&0&0\cr0&0&0&0&0&0\cr0&0&0&0&0&0\cr\end{matrix}\right)\deuf$$\normalsize
where
$$\widetilde{k_2}= \dx(k_1)+\frac{5}{2}L_1, \quad  \widetilde{k_3}=\dy(k_1)+\frac{5}{2}L_2,\quad\widetilde{k_4}=\partial^2_x(k_1)+4\dx(L_1)+\frac{v_4}{2}L_1-\frac{3}{2}v_5L_2,$$
$$\widetilde{k_5}=\dx\dy(k_1)-2\dx(L_2)+2\dy(L_1)+\frac{v_3}{2}L_1-\frac{v_4}{2}L_2\quad\widetilde{k_6}=\partial^2_y(k_1)+4\dy(L_2)+\frac{3}{2}v_2L_1-\frac{v_3}{2}L_2,$$
and where~$<v_1>_i$ belongs to the differential ideal generated by
$v_1$ so that if~$v_1=0$, then ~$<v_1>_i=0$. Moreover, the
curvature is an $\ol^*$-invariant of the web.
\end{theo}

This result is quite similar to the case of $4$-webs.
Nevertheless, the existing difference makes clearer the fact that
Poincar\'e theorem is not valid for $5$-webs. Instead, we see in
the following theorem that $v_1$ is an obstruction for a maximal
rank web to be linearizable:\begin{theo}\label{vex} Let~$\W(5)$ be
a~$5$-web of maximal rank. Then $\W(5)$ is linearizable if and
only if~$deg\;P_{\W(5)}\leq 3$.
\end{theo}
\begin{proof}
If~$deg\left(P_{\W(5)}\right)\leq3$, then~$v_1=0$ and the
curvature is
$$K=\left(\begin{matrix}k_1&\widetilde{k_2}&\widetilde{k_3}&\widetilde{k_4}&\widetilde{k_5}&\widetilde{k_6}\cr0&0&0&0&0&0\cr0&0&0&0&0&0\cr0&0&0&0&0&0\cr0&0&0&0&0&0\cr0&0&0&0&0&0\cr\end{matrix}\right)\deuf$$
Since the web is of maximal rank, the curvature is equal to zero.
So~$k_1=0$ and then~$L_1$ and~$L_2$ too. The properties \ref{pw}
gives us the conclusion. The other implication is also a direct
consequence of the same properties \ref{pw}.\end{proof} This
theorem was first proved in \cite{H-3}, but we give here a new proof only using the connection.

\subsection{Revisiting exceptional~$5$-webs} An exceptional web is
a web of maximal rank which is not algebrizable, or with the
theorem of Lie-Darboux-Griffiths, it is the same to say that it is
not linearizable. The Poincar\'e theorem for~$4$-webs says that this
configuration is only possible for~$d\geq5$.

Our expression of the connection offers a new approach of
exceptional webs. We do not need to exhibit abelian relations to
determine wether a web is exceptional or not. Moreover, all
exceptional~$5$-webs can be described thanks to a differential
system we give now:

\begin{prop}\label{5w} Let~$\W(5)$ be a planar~$5$-web. We have the following equivalence\dept
\begin{enumerate}
    \item[i)]$\W(5)$ is exceptional\pvir
    \item[ii)] The following explicit conditions are satisfied\dept
   ~$$\la\begin{matrix}k_1=k_2=k_3=k_4=k_5=k_6=0\cr v_1\neq 0\,\in\ol\end{matrix}\right..$$
\end{enumerate}
\end{prop}
This system can be seen as systems in the coefficients~$a_i$ of
the differential equation which presents the web. Hence, the study
of exceptional~$5$-webs could be lead theoretically through the
study of this system. But clearly, its complexity draws the limits
of such an approach.

\section{The trace formula}
The trace formula links the trace of the curvature associated with
a $d$-web~$\W(d)$, namely~$k_1$, to the Blaschke
curvatures~$d\gamma_k$ of the~$3$-webs extracted
from~$\W(d)$, $(^d_3)$ in number:
\begin{theo}[Trace formula]\label{trace} Let~$\W(d)$ a planar~$d$-web. With the previous notations,
~$$\mbox{\upshape tr \itshape} K=k_1=\sum_{k=1}^{(^d_3)}
d\gamma_k,$$
\end{theo}
This formula was first demonstrated for~$4,5$ and~$6$ webs in the
author thesis \cite{R2}. In this continuity, we will give here a
full proof of this theorem in subsection \ref{ptf}.

The formula has various interpretations and consequences. Still in
the author thesis, we made a construction of a poly-hexagon that
generalizes the one of Thomsen for $3$-webs, based on this trace
formula. It gives also a simple but very useful criteria for
searching exceptional webs.

One must notice that in 1938, Pantazi (\cite{Pa}) gave a
construction to determine the maximal rank webs which lead him to
the introduction of~$\pi_d$ expressions whose annulation gives the
conditions for a web to be of maximal rank. Mihaileanu
(\cite{Mih}), following Pantazi, identified one of this
coefficient to be the sum of the Blaschke curvatures of
extracted~$3$-webs. Our results were stated independently from
them, since they were ignored until Luc Pirio digs them in his
thesis. These results offer obvious similarity with the previous
statements, even if the construction is not the same and Pantazi
and Mihaileanu did not give a truly proof of there results. But
the links between both approaches can be found in \cite{H2R},
unifying the trace formulas. Note also that a general proof of the
trace formula which use another, but
equivalent, formalism will be given in the same paper.

Here, the trace formula will always refer to the relation
satisfied by the trace of the curvature introduced here.

$~$

In order to produce our proof for a $d$-web, we must first exhibit
the relations that link, say the coefficients of the connection
associated with the extracted~$3$-webs and the coefficients of
the~$d$-web itself. We then could calculate the sum of the
Blaschke curvatures in terms of the coefficients of the connection
of the~$d$-web. This is our first step. Then we need to
give a general expression of the trace~$k_1$ of the curvature, and
compare the two expressions which will be the second step.

Before we prove the trace formula, let us make a remark.
We let~$\W(d)$ be a~$d$-web and its
presentation~$F(x,y,y')=0$ . Out of the singular locus, the
slopes of the web will be denoted by $p_i\in\ol$, for~$1\leq i\leq
d$. Given such a~$d$-web, we want to know wether an abelian relation comes from extracted webs of $\W(d)$.

Let~$p_k$ be one of the slope of the web, and~$\W_k(d-1)$ be the
extracted~$d-1$ web of~$\W(d)$ obtained by forgetting the
slope~$p_k$. We want to know the $\ol^*$-invariants associated
with this web thanks to those of the~$d$-web.

On can show (\cite{H-1}) that the isomorphism~$T$
between~$\mathcal{A}(d)$ and $\mathfrak{a}_F$ associates to an
abelian relation $(g_i(F_i))_{\la 1\leq i\leq d\ra}$ the~$1$-form
on~$S$:$$T( (g_i(F_i))_{1\leq i\leq d})=r(x,y,p)\fs \mbox{ where }
r=F.\left(\sum_{i=1}^d\frac{g_i(F_i)\dy(F_i)}{p-p_i}\right).$$
\begin{prop}\label{compl}An abelian relation of~$\W(d)$ is an abelian relation of the web~$\W_k(d-1)$ if and only if
$r(x,y,p_j)=0$ \emph{i.e.}~$r=(p-p_k).r_k$ where~$r_k$ is a
polynomial in~$p$ of degree~$d-1$.
\end{prop}
\begin{proof} An abelian relation~$(g_i(F_i))_{\la 1\leq i\leq d\ra}$ of~$\W(d)$ is an abelian relation of
$\W_k(d-1)$ if and only if~$g_k(F_k)=0$. Via the isomorphism~$T$
this is a necessary and sufficient condition for~$r$ to
admit~$p_k$ as a root.
\end{proof}
\subsection{Extracted webs}
Let~$p_k$ be one of the slope of the web, and~$\W_k(d-1)$ be the
extracted~$d-1$ web of~$\W(d)$ obtained by forgetting the slope
$p_k$. It is presented by~$F_k(x,y,p)=0$ with the relation
$$F(x,y,p)=\lp p-p_k(x,y)\rp F_k(x,y,p).$$

We want to compute the coefficients of the system $\Mr(d-1)$ associated with the web~$\W_k(d-1)$ with the ones of $\Mr(d)$, associated with $\W(d)$. Let~$S$ be the surface of~$\cl^3$ defined by
$S=~\la~F(x,y,p)~=~0~\ra$ and~$S_k$ the surface defined by
$S_k=\la F_k(x,y,p)=0\ra$ and let $r_k(x,y,p)$ be a polynomial in~$p$ with coefficients in~$\ol$
of degree~$d-4$. The canonical monomorphism ~$i:S_k\rightarrow S$
induces the
morphism~$$\begin{matrix}i^*:&\Omega^1_S=\Omega^1_{\cl^3}/\left(
dF,F\Omega^1_{\cl^3}
\right)&\rightarrow&\Omega^1_{S_k}=\Omega^1_{\cl^3}/\left(
dF_k,F_k\Omega^1_{\cl^3}\right)\cr & (p-p_k).r_k(x,y,p)\fs&
&\omega_k=r_k(x,y,p)\fsk\end{matrix}$$

Indeed, we have
$\parp(F)=\parp(F_k)(p-p_k)+F_k$ where~$F_k=0$ on~$S_k$
and~$$i^*(\omega)=(p-p_k).r_k(x,y,p)\circ
i\frac{dy-p\;dx}{\parp(F_k)(p-p_k)\circ\; i+F_k\circ i}$$
by definition of pull back.
By the same way, we show that for a polynomial~$t$  in~$\ol[p]$,
we have
$$i^*(t\frac{dx\wedge dy}{\parp(F)})=t\frac{dx\wedge
dy}{(p-p_k)\parp(F_k)}.$$

As reminded in section $2$, given a~$1$-form~$\sigma$ on~$S$, there exists a
polynomial $t_\sigma$ of degree~$d-1$ in~$\ol[p]$ such that the differential
can be written
$\displaystyle{d\sigma=t_\sigma\frac{dx\wedge dy}{\parp(F)}}$.
We thus have the following proposition :
\begin{prop}Let $\omega_k$ be a $1$-form on $S_k$ and $t_{\omega_k}$ the polynomial of degree $d-1$ such that $\displaystyle{d\omega_k=t_{\omega_k}\frac{dx\wedge dy}{(\parp(F_k)}}$. Let $\omega=(p-p_k).\omega_k$ be the corresponding $1$-form on $S$ by $i^*$. Then we have $$d\omega=(p-p_k)\cdot t_{\omega_k}\frac{dx\wedge
dy}{\parp(F)}.$$
\end{prop}
\begin{proof} Indeed, we have the following equalities:
$\displaystyle{d(i^*(\omega))=d\omega_k=t_{\omega_k}\frac{dx\wedge dy}{(\parp(F_k)}}$
and $\displaystyle{i^*(d\omega)=i^*(t_\omega\frac{dx\wedge
dy}{\parp(F)})=t_\omega\frac{dx\wedge dy}{(p-p_k)\parp(F_k)}}.$ Since the differential and the pull back commutes
$d(i^*(\omega))=i^*(d\omega)$, so~$t_\omega=t_{\omega_k}.(p-p_k).$ \end{proof}
 Notice that this is another proof of proposition
\ref{compl}, by taking $t_{\omega_k}=0$.

If we let~$r_k(x,y,p)=b_3p^{d-4}+\ldots+b_{d-1}$, the differential on~$S_k$ of $\omega_k=r_k\cdot\fs$ induces a system denoted as follow, where the $t_i$ are the coefficients of $t_{\omega_k}$:
\small$$S(A^k_{ij}): \left\{\begin{matrix}\dx (b_{d-1}) &+&
A^k_{1,1}\cdot b_3 + \cdots + A^k_{1,d-3}\cdot
b_{d-1}&=&t_{d-1}\cr \dx (b_{d-2})+\dy(b_{d-1})&+& A^k_{2,1}\cdot
b_3 + \cdots + A^k_{2,d-3}\cdot b_{d-1}&=&t_{d-2}\cr
  & & & \vdots& \cr
\dx (b_{3})+\dy(b_{4})&+&  A^k_{d-3,1}\cdot b_3 + \cdots +
A^k_{d-3,d-3}\cdot b_{d-1}&=&t_3\cr \dy(b_3)&+&A^k_{d-2,1}\cdot
b_3 + \cdots + A^k_{d-2,d-3}\cdot
b_{d-1}&=&t_2\cr\end{matrix}\right.$$ \normalsize\indent
So the corresponding system for the differential on $S$ is then \small$$\left\{\begin{matrix}\dx (-p_k.b_{d-1}) +
A_{1,1}\cdot b_3 + \cdots + A_{1,d-2}\cdot
(-p_k.b_{d-1})&=&-p_k.t_{d-1}\cr \dx
(b_{d-1}-p_k.b_{d-2})+\dy(-p_k.b_{d-1})+   A_{2,1}\cdot b_3 +
\cdots + A_{2,d-2}\cdot (-p_k.b_{d-1})&=&t_{d-1}-p_k.t_{d-2}\cr
   & \vdots& \cr
\dx (b_{3})+\dy(b_{4}-b_3.p_k)+  A_{d-2,1}\cdot b_3 + \cdots +
A_{d-2,d-2}\cdot (-p_k.b_{d-1})&=&t_3-p_k.t_2\cr
\dy(b_3)+A_{d-1,1}\cdot b_3 + \cdots + A_{d-1,d-2}\cdot
(-p_k.b_{d-1})&=&t_2\cr\end{matrix}\right.$$ \normalsize

By a combination of the preceding systems, but omitting the
explicit expressions, we give the following claim:
\begin{theo}\label{extr}With the previous development, we get an
explicit writing of the coefficients~$A^k_{ij}$ of the extracted
web thanks to the~$A_{ij}$ and the forgotten slope.
\end{theo}

\begin{ex}\normalfont For a~$4$-web, we have the following relations:
$$\la\begin{matrix}A^k_2=A_2-v_2-v_1p_k \cr
A_1^k=-\dy(p_k)+A_1-v_3-v_2p_k-v_1p_k^2\cr
-V(p_k)=\dx(p_k)+p_k\dy(p_k)\end{matrix}\right..$$

In the linear case, where~$A_{11}=A_{32}=0$,~$A_{12}=A_{21}=A_1$,
and $A_{22}=A_{31}=A_2$, we get~$A_2=A_2^k$ and
$A_1=A_1^k-\dy(p_k),$ which means that the trace formula is
directly checked, since the Blaschke curvatures are the
$2$-forms~$(\dx(A_2^k)-\dy(A_1^k))\deuf$.\end{ex}

\subsection{Proof of the trace formula}\label{ptf}
The expression of the coefficients ~$(A^k_1,A^k_2)$ of all the extracted~$3$-webs of a
$d$-web allows us after a long computation to give the expression
of the sum of the Blaschke curvatures in terms of the coefficients
$A_{ij}$ associated with the~$d$-web. This can be done using heavily
the Newton relations between the roots and the coefficients of the
equation presenting the web.

The result of this sum is a quite simple expression, which needs to be the
trace of the curvature matrix in order to prove the trace formula.
\begin{prop}The trace of the connection, in
a suitable basis is given by:
$$
\mbox{\upshape tr \itshape}\gamma=\sum_{q=1}^{d-2}
A_{d-q-1,q}dx+A_{d-q,q}dy+
\sum_{q=2}^{d-2}A_{d-q-1,q}dx+A_{d-q+1,q-1}dy$$
\end{prop}
Again the proof is a computation. Since the
differential of this trace is the trace of the curvature
matrix~$\mbox{tr }K$, it does not depends on the choice of the
basis. We then compute more easily this trace by choosing a
suitable basis for this.

Moreover, we do not need to have a general expression of the basis,
since the construction of the adapted basis constrain the vectors
of the basis to have a particular form which is sufficient to
compute the trace. This trace is, as expected, the sum of the Blaschke curvatures computed before.

We give now an interpretation of the trace formula in terms of
determinant:
\begin{theo}\label{det} Let~$\W(d)$ be a planar~$d$-web
 and let~$(L_k,\pq_k)$ be the line bundles associated with extracted~$3$-webs of~$\W(d)$ for~$1\leq k\leq (^d_3)$.
We have an isomorphism of line bundles with connection:
$$(\det\,E,\det\,\pq)\;\widetilde{=}\;(\bigotimes_{k=1}^{(^d_3)}\,L_k,\bigotimes_{k=1}^{(^d_3)}\,\pq_k),$$
\end{theo}
\begin{proof}The isomorphism denoted~$\tau$ is defined by its action on the basis.
To the basis~$e_1\wedge\ldots\wedge e_{\pi_d}$ of~$\det\,E$ which
is the wedge product of the vector of the basis of~$E$, we
associate the basis~$\bigotimes_ke^k$ of
$\bigotimes_{k=1}^{(^d_3)}\,L_k$ where~$e^k$ is a basis of~$L_k$.
The isomorphism commutes with the connections since the equality
$$\bigotimes_k\gamma_k\circ \tau\left(\bigwedge_k
e_k\right)=\tau\circ \det\,\pq\left(\bigwedge_k e_k\right)$$ is a
consequence, regarding matrices, of the trace formula
$$\mbox{\upshape tr \itshape}\gamma\otimes\left(\bigotimes_ke^k\right)=\left(\sum_{k=1}^{(^d_3)}
\gamma_k\right)\otimes\left(\bigotimes_ke^k\right).$$
\end{proof}

\section{Determination of the rank}
\subsection{The main result}

The  determination of the rank of a $d$-web was not effective since Blaschke
introduced the subject, except for $d=3$ as we have seen. In fact, the main way used before to get the rank was
to compute all the abelian relations.

The following theorem gives, not only a determination of the rank,
but also an explicit locally free~$\ol$-module of finite type
whose rank is the rank of the web. It embodies the abelian
relations, since it is generated by them. The existence of such a
$\ol$-modules is theoretically given, but the fact that it can be embodied
explicitly is not clear in general. We give here a complete proof of this result first announced in
our note \cite{R1}:
\begin{theo}[Determination of the rank]\label{rk} Let~$\W(d)$ a non singular planar~$d$-web. There
exists a $\cl$-vector bundle~$\overline{K}$ of rank~$\mbox{\upshape rk
\itshape}\W(d)$, which is the kernel of an explicit endomorphism
of~$\ol^{\pi_d}$ such that
$$\overline{K}=\ol\otimes_\cl \mbox{\upshape Ker \itshape}\pq.$$ \indent So if
~$(k_{m\ell})$ denotes the matrix of this endomorphism, the rank
of the web is given by\upshape:
\itshape$$\mbox{\upshape rk \itshape}\W(d)=\mbox{\upshape corank
\itshape}(k_{m\ell}).$$
\end{theo}
Before we prove this theorem, we give the construction of the
matrix~$(k_{ml})$. The horizontal sections $f=\,^t(f_1, f_2, \ldots
,f_{\pi_d})\in E^\pq:=\mbox{Ker }\pq$ of~$\pq$ are identified to
the abelian relations of the web by construction of the connection,
and satisfy the differential system~$df+\gamma f=0$ where~$\gamma=\gamma_xdx+\gamma_ydy$
is the connection matrix in a suitable basis. The integrability
condition is then given by the \emph{only} one relation~$$(1)\quad k\cdot
f=k_1f_1+k_2f_2+\cdots+k_{\pi_d}f_{\pi_d}=0$$ where the
$k_i$ are the coefficients of the curvature matrix. We consider then the~$\pi_d$ equations obtained from the
derivation until the order~$d-3$ of~$(1)$ where we substitute the
derivative of~$f$ thanks to~$df=-\gamma f$. We get then a square
matrix~$(k_{m\ell})$ of order~$\pi_d$ whose first line is the
first line of the curvature matrix:

$$\lp\begin{matrix}k_1&k_2&\ldots & k_{\pi_d} \cr k_{21} &k_{22} &\ldots
& k_{2,\pi_d} \cr \vdots & \vdots&  &\vdots \cr k_{\pi_d,1} &k_{\pi_d,2}&
\ldots  &k_{\pi_d,\pi_d} \cr\end{matrix}\rp$$

\begin{proof}
We will largely use the fact that~$E^\pq$ is a local system. Let
$\overline{K}:=\mbox{Ker}\,(k_{m\ell})$ the~$\ol$-module of finite
type  defined by our matrix. We have the inclusion
$$\ol\otimes_\cl E^\pq\subseteq \overline{K}\qquad (2)$$ by
construction. The converse will be true by Nakayama's lemma if we
show the following
$$\overline{K}=\ol\otimes_\cl
E^\pq+\mathfrak{m}\cdot\overline{K}\qquad (3)$$ where
~$\mathfrak{m}$ is the maximal ideal of~$\ol$. Let
$(g_1,\ldots,g_r)$ a system of generators of~$\overline{K}$ such
that the~$g_i(0)$ are linearly independent and whose existence is
again given by Nakayama's lemma.

Let~$g\in\overline{K}$. One can write that~$g=\sum_{1\leq i\leq
r}\lambda_ig_i$ with  $\lambda_i\in\ol.$ If~$g(0)=0$, then
$g(0)=\sum_{1\leq i\leq  r}\lambda_i(0)g_i(0)=0$ and
so~$\lambda_i(0)=0$ since the~$g_i(0)$ are linearly independent.
This gives~$g\in\mathfrak{m}\overline{K}$ and the needed equality.

If~$g(0)\neq0$, we will construct an analytic function ~$f$
in~$E^\pq$ such that~$f(0)=g(0)$, where~$g(0)$ stands for the
initial conditions in the classical Cauchy Theorem.

Indeed, Cauchy theorem in one variable gives us a unique function
$\sigma(y)$ such that~$\sigma(0)=g(0)$ and
$\dy(\sigma)+\gamma_y\sigma=0$. Again, Cauchy theorem gives us a
unique function~$f(x,y)$ such that~$f(0,y)=\sigma(y)$ and
$\dx(f)+\gamma_xf=0$. We have the following two equalities
:~$$\la\begin{matrix}\dy(f(0,y))+\gamma_yf(0,y)&=&0&&&&(4)\cr\dx(f(x,y))+\gamma_xf(x,y)&=&0&&&&(5)\end{matrix}\right.$$
Let~$\tau=\dy(f)+\gamma_yf$. The equation~$(4)$ gives us
that~$\tau(0,y)=0$, and using~$(5)$, we
have~$\dx(\tau)+\gamma_x\tau=k\cdot f.$

Since~$g\in\overline{K}$ and~$f(0)=g(0)$, we have~$k(0)\cdot
f(0)=0$, and all the successive derivatives of~$k\cdot f$ taken at
$0$ are zero, again by~$(4)$, ~$(5)$ and the
equality~$\overline{K}(0)\cdot g(0)=0$. So we have~$k\cdot f=0,$
and then,~$\dx(\tau)+\gamma_x\tau=0$ and $\tau(0,y)=0.$ The
uniqueness theorem of Cauchy gives again~$\tau=\dy(f)+\gamma_yf=0$
and so~$f$ belongs to~$E^\pq$ since we
have~$$\la\begin{matrix}\dy(f)+\gamma_yf&=&0\cr\dx(f)+\gamma_xf&=&0\end{matrix}\right..$$
The inclusion~$(2)$ gives~$g-f\in\overline{K}$ and~$(g-f)(0)=0$.
But we have seen that if~$h\in\overline{K}$ is such that~$h(0)=0$,
the system of generators of~$\overline{K}$  chosen allows us to
show that~$h\in\mathfrak{m}\overline{K}$.
So~$g-f\in\mathfrak{m}\overline{K}$, which shows the needed
equality.
\end{proof}As a corollary, we precise the fact that a web is, in general, of rank equal to~$0$:
\begin{cor}A~$d$-web is of rank greater or equal to~$1$ if and only if~$det(k_{ml})=0$
\end{cor}

\begin{ex}\normalfont For~$d=4$, we have the explicit expression of the coefficients of the matrix
$(k_{ml})$:
$$k_{21}=\dx(k_1)-A_1k_1+k_2\quad k_{22}=\dx(k_2)-(\dy(v_4)+v_4v_2)k_1-(A_1-v_3)k_2-v_4k_3$$
$$k_{23}=\dx(k_3)-(v_1v_4-\kappa_2)k_1-A_1k_3\quad k_{31}=\dy(k_1)-(A_2-v_2)k_1+k_3$$
$$k_{32}=\dy(k_2)-(v_1v_4+\kappa_1)k_1-(A_2-v_2)k_2 \quad k_{33}=\dy(k_3)-(v_1v_3-\dx(v_1))k_1+v_1k_2-A_2k_3$$
\end{ex}
\subsection{Linear case for~$4$-webs}
\begin{prop} The rank of a linear~$4$-web is not~$2$.\end{prop} Indeed, let~$\mathcal{L}(4)$ be a linear~$4$-web. The
Lie-Darboux-Griffiths says that if such a web admits a
complete abelian relation, \emph{i.e.}~$(g_i(F_i))_{1\leq i\leq
4}$ where~$g_i$ is non zero for all~$i$, then the web is of
maximal rank.

If~$\mathcal{L}(4)$ is of rank 2 and so does not have any complete
abelian relations, this web admits two independent abelian relations
coming from (different) extracted~$3$-webs.In fact, there exists a linear combinaison of these two
abelian relations which is then complete, in
contradiction with the hypothesis. Thus, a linear~$4$ web can not be of rank~$2$.

We can also prove this using our preceding results: the
cancelation of all the minors of order~$1$ and~$2$ implies that
$k_1=0$ which is impossible since this implies for a linear web to
be of maximal rank.

\subsection{Bol's theorem for~$4$-webs}
As a nice application of the trace formula and the theorem of
determination of the rank, we can prove the Bol's theorem for
4-webs. But before, we have this rigidity type proposition:
\begin{prop} A~$4$-web of rank at least~$2$ with~$k_1=0$ is of
maximal rank.
\end{prop}
\begin{proof}
A computation show that if~$k_1=0$ and all the
minors of order~$1$ and~$2$ of~$k_{ml}$ are equal to zero, then the matrix of curvature
is zero. Hence the web is of maximal rank.
\end{proof}
\begin{theo}
Let~$\mathcal{H}(4)$ be a \emph{hexagonal}~$4$-web (\emph{i.e.} the
Blaschke curvatures of all extracted~$3$-webs are zero). Then
$\mathcal{H}(4)$ is of maximal rank.
\end{theo}
\begin{proof}
Since the web is hexagonal, the trace formula shows that~$k_1=0$
and there exists at least two abelian relations of~$3$-webs
linearly independent. So the rank is at least~$2$ and the
preceding proposition applies. \end{proof}

\section{Example}

Let~$\W(4)$ be the~$4$-web presented by
${p}^{4}+{y}^{2}{p}^{2}-yp=0$. Its discriminant is
$\Delta=-{y}^{4} \left( 27+4{y}^{4} \right)$ and
$$P_{\mathcal{W}(4)}=-12{\frac {{p}^{3}}{27+4{y}^{4}}}+{\frac {
\left( 9+4{y}^{4}
 \right) {p}^{2}}{y \left( 27+4{y}^{4} \right) }}-8{\frac {{y}^{2}
p}{27+4{y}^{4}}}
$$ The fundamental forms is $\alpha=-2{\frac { \left( 9+4{y}^{4} \right)}{y \left( 27+4{y}^
{4} \right) }}dy.$ We can compute the connexion matrix:
$$\gamma=\left( \begin {array}{ccc} -{\frac { \left( 9+4{y}^{4} \right) { dy}}{y \left( 27+4{y}^{4} \right) }}&-16{\frac {y \left( -27+4
{y}^{4} \right) { dy}}{ \left( 27+4{y}^{4} \right) ^{2}}}&96{
\frac {{y}^{2}{ dy}}{ \left( 27+4{y}^{4} \right) ^{2}}}
\\\noalign{\medskip}-{ dx}&-8{\frac {{y}^{2}{ dx}}{27+4{y}^{
4}}}-{\frac { \left( 9+4{y}^{4} \right) { dy}}{y \left( 27+4{y}
^{4} \right) }}&-12{\frac {{ dy}}{27+4{y}^{4}}}
\\\noalign{\medskip}-{ dy}&0&-2{\frac { \left( 9+4{y}^{4}
 \right) { dy}}{y \left( 27+4{y}^{4} \right) }}\end {array}
 \right)
$$
and the curvature matrix:
$$K=\left( \begin {array}{ccc} -16{\frac {y \left( -27+4{y}^{4} \right) }{ \left( 27+4{y}^{4} \right) ^{2}}}&-128{\frac {{y}^{3}
 \left( -27+4{y}^{4} \right) }{ \left( 27+4{y}^{4} \right) ^{3}}}&0
\\\noalign{\medskip}0&0&0\\\noalign{\medskip}0&0&0\end {array}
 \right)~$$
The rank of this web is~$2$ since the matrix~$(k_{m\ell})$ is the
following:
$$(k_{m\ell})= \left(
\begin {array}{ccc} -16{\frac {y \left( -27+4{y}^{4} \right) }{
\left( 27+4{y}^{4} \right) ^{2}}}&-128{\frac {{y}^{3}
 \left( -27+4{y}^{4} \right) }{ \left( 27+4{y}^{4} \right) ^{3}}}&0
\\\noalign{\medskip}-128{\frac {{y}^{3} \left( -27+4{y}^{4}
 \right) }{ \left( 27+4{y}^{4} \right) ^{3}}}&-1024{\frac {{y}^{5}
 \left( -27+4{y}^{4} \right) }{ \left( 27+4{y}^{4} \right) ^{4}}}&0
\\\noalign{\medskip}64{\frac {243-306{y}^{4}+8{y}^{8}}{ \left(
27+4{y}^{4} \right) ^{3}}}&512{\frac {{y}^{2} \left( 243-306{y}^
{4}+8{y}^{8} \right) }{ \left( 27+4{y}^{4} \right) ^{4}}}&0
\end {array} \right)~$$ whose determinant and order~$2$ minors are all equal to zero.


\begin{thebibliography}{aaaaa}
\bibitem[B]{B}W. Blaschke, Einfuhrung in die Geometrie der Waben, Birkhauser, Basel, 1955.
\bibitem[B-B]{BB}W. Blaschke und G. Bol, Geometrie der Gewebe, Springer, Berlin, 1938.
\bibitem[BC3G]{BC3G}R.L. Bryant, S.S. Chern, R.B. Gardner, H.L. Goldschmidt and  P.A. Griffiths, Exterior differential systems, Springer, Berlin, 1991.
\bibitem[C]{C}S.S. Chern, Web Geometry, Bull. Amer. Math. Soc.  \bfseries6 \normalfont (1982), 1-8.
\bibitem[C-G]{CG}S.S. Chern and  P.A. Griffiths, Abel's Theorem and Webs, Jahresber. Deutsch. Math.-Verein.  \bfseries80 \normalfont (1978), 13-110 and Corrections and Addenda to Our Paper: Abel's Theorem and Webs,  Jahresber. Deutsch. Math.-Verein.  \bfseries83 \normalfont (1981),78-83.
\bibitem[C-L]{CL}V. Cavalier et D. Lehmann, Introduction \`a une \'etude globale des tissus sur une surface holomorphe, To appear.
\bibitem[G-02]{G2}P.A. Griffiths, The legacy of Abel in algebraic geometry, \textit{in} Laudal, Olav Arnfinn (ed.) et al., \textit{The legacy of Niels Henrik Abel. Papers from the Abel bicentennial conference, University of Oslo, Oslo, Norway, June 3-8, 2002}. Springer, Berlin (2004), 179-205.
\bibitem[G-76]{G}P.A. Griffiths, Variations on a Theorem of Abel, Invent. Math. \bfseries35 \normalfont (1976), 321-390.
\bibitem[Go]{Go}V.V. Goldberg, 4-webs in the plane and their linearizability  Acta. Appl. Math. \textbf{80} (2004), 35-55.
\bibitem[Go-L]{GL}V.V. Goldberg and V.V. Lychagin, On linearisation of planar three-webs and Blaschke's conjecture, C.R. Acad. Sci. Paris, Ser. I \bfseries341  \normalfont (2005), 169-173.
\bibitem[H2R]{H2R}A. H\'enaut, O. Ripoll et G. Robert, Formule de la trace pour la connexion d'un tissu du plan, in preparation.
\bibitem[H-07]{H-2}A. H\'enaut, Planar web geometry through abelian relations and singularities, \emph{in} Inspired by Chern, A memorial volume in honor of a great mathematician, World Scientific, Sci. Publishing co., River Edge, NJ, 2006.
\bibitem[H-04]{H-1} A. H\'enaut, On planar web geometry through abelian relations and connections, Ann. of Math. \textbf{159} (2004), 425-445.
\bibitem[H-94]{H-3}  A. H\'enaut, Caract\'erisation des tissus de~$\cl^2$ dont le rang est maximal et qui sont lin\'earisables, Compositio Math. \bfseries94 \normalfont (1994), 247-268.
\bibitem[M2P]{PPM}  D. Mar\'in, J.V. Pereira and L. Pirio, On planar webs with infinitesimal automorphisms, \emph{in} Inspired by Chern, A memorial volume in honor of a great mathematician, World Scientific, Sci. Publishing co., River Edge, NJ, 2006.
\bibitem[Mih]{Mih}  M.N. Mihaileanu, Sur les tissus plans de premi\`ere esp\`ece, Bull. Math. Soc. Roum. Sci.  \textbf{43} (1941), 23-26.
\bibitem[N]{N}  I. Nakai, Curvature of curvilinear 4-webs and pencils of one forms: Variation on a theorem of Poincar\'e, Mayrhofer and Reidemeister, Comment. Math. Helv.  \bfseries73 \normalfont (1998), 177-205.
\bibitem[P]{P}  L. Pirio, \'Equations fonctionnelles ab\'eliennes et g\'eom\'etrie des tissus, Th\`ese de doctorat, Universit\'e Paris VI, d\'ecembre 2004.
\bibitem[P-T]{PT}   L. Pirio et J.-M. Tr\'epreau, Tissus plans exceptionnels et fonction th\^eta, Annales de l'Institut Fourier, \textbf{55} no. 7 (2005), 2209-2237.
\bibitem[Pa]{Pa}  A. Pantazi, Sur la d\'etermination du rang d'un tissu plan, C.R. Acad. Sci. Roumanie \bfseries4  \normalfont (1938), 108-111.
\bibitem[R1]{R1}  O. Ripoll, D\'etermination du rang des tissus du plan et autres invariants g\'eom\'etriques, C.R. Acad. Sci. Paris, Ser. I \bfseries341  \normalfont (2005), 247-252.
\bibitem[R2]{R2}  O. Ripoll, G\'eom\'etrie des tissus du plan et \'equations diff\'erentielles, Th\`ese de doctorat, Universit\'e Bordeaux I, d\'ecembre 2005, disponible sur http://tel.archives-ouvertes.fr/tel-00011928.
\bibitem[R-S-1]{RS}  O. Ripoll et J. Sebag, Solutions singuli\`eres des tissus polynomiaux du plan, J. of Algebra \textbf{310}, (2007), 351-370.
\bibitem[R-S-2]{RS2}  O. Ripoll et J. Sebag, Tissus du plan et polynômes de Darboux, (submitted).
\bibitem[W]{W}  J. Grifone and  \'E. Salem (Eds), Web Theory and Related Topics, World Scientific, Sci. Publishing co., River Edge, NJ, 2001.
\end{thebibliography}
\end{document}